\documentclass[reqno]{amsart}

\usepackage{enumerate}
\usepackage{amstext,amsmath,amsthm,amsfonts,amssymb,amscd}
\usepackage{latexsym,mathrsfs,eufrak,dsfont,euscript}

\usepackage{hyperref}

\usepackage{tikz}
\usepackage{pgfplots}
\usetikzlibrary{calc,through,backgrounds}
\usetikzlibrary{calendar,snakes}
\usetikzlibrary{positioning}

\newtheorem{theorem}{Theorem}

\theoremstyle{definition}

\theoremstyle{remark}

\numberwithin{equation}{section}
\newcommand{\norm}[1]{\left\Vert#1\right\Vert}
\begin{document}
\title[]{A note on Besov regularity for parabolic initial boundary value problems}

\author[]{Hugo Aimar}
\email{haimar@santafe-conicet.gov.ar}

\author[]{Ivana G\'{o}mez}
\email{ivanagomez@santafe-conicet.gov.ar}
\thanks{The research was supported  by CONICET, ANPCyT, and UNL}



\begin{abstract}
In this note we consider the initial boundary value problem for the heat equation on cylinders based on Lipschitz domains with Besov data. We obtain a regularity exponent for the solution that improves the rate of convergence of nonlinear approximation methods.
\end{abstract}

\maketitle

\section{Introduction}
As in the elliptic results by Dahlke and DeVore in \cite{DaDeV97}, based in \cite{DJP92}, the improvement of the regularity exponent in Besov norms becomes a tool to measure the rate of convergence of nonlinear approximation methods. A parabolic Besov regularity improvement for temperatures, i.e., solutions of $\tfrac{\partial u}{\partial t}=\Delta u$, follows from the results obtained by the authors in \cite{AGconstapp}, \cite{AGIjfaa} and \cite{AGIjfa}. Nevertheless to apply those results, a starting Besov regularity for the temperature is required. In the elliptic case this initial regularity for the harmonic function is proved to be attained if we solve the Dirichlet problem with Besov boundary data, see the results of Jerison and Kenig in \cite{JeKe95}. Its parabolic counterpart is due to Jakab and Mitrea, see \cite{JaMit06}.

\medskip
For a given $0<\varepsilon\leq 1$ we shall write $\mathcal{R}_{\varepsilon}$ to denote the set of those points $(a,b)\in [0,1]^2$ in the plane satisfying any one of the following three conditions
\begin{align*}
\frac{1-\varepsilon}{2}<b<\frac{1+\varepsilon}{2} &\textrm{ and } 0<a<1;\\
\frac{1+\varepsilon}{2}\leq b<1 &\textrm{  and } 2b-1-\varepsilon<a<1;\\
0<b\leq\frac{1-\varepsilon}{2} &\textrm{  and } 0<a<2b+\varepsilon.
\end{align*}

The following figure contains a picture of $\mathcal{R}_{\varepsilon}$ for $\varepsilon=\tfrac{1}{4}$.
\newpage

\begin{figure}[h]
\begin{center}
\begin{tikzpicture}[scale=3]
\draw[line width=.5pt, dotted, fill=lightgray!40]
(0,0)--(1/4,0)--(1,1/4)--(1,1)--(3/4,1)--(0,5/8)--(0,0);
\draw[->] (-0.01,0)  -- (1.2,0) node[right] {\scalebox{1}{$\displaystyle a$}};
\draw[->] (0,-0.01) -- (0,1.2) node[left] {\scalebox{1}{$b$}};
\draw[shift={(0,1)}] (.5pt,0pt) -- (-.5pt,0pt) node[left] {\scalebox{.7}{$1$}};
\draw[shift={(1,0)}] (0pt,-0.5pt)--(0pt,0.5pt) node[below] {\scalebox{.7}{$1$}};
\draw[line width=.7pt, snake=brace] (0,0.05)--(1/4,0.05);
    \node[line width=.5pt] at (1/8,.14) {{\scalebox{1}{$\varepsilon$}}};
\draw[line width=.7pt, snake=brace] (-.05,1/2)--(-0.05,5/8);
    \node[line width=.5pt] at (-0.14,0.55) {\scalebox{1}{$\tfrac{\varepsilon}{2}$}};
\end{tikzpicture}
\end{center}
\caption{The region $\mathcal{R}_{\varepsilon}$ for the parameters $(s,\tfrac{1}{p})$ of the regula\-ri\-ty of data.}\label{fig:regiondata}
\end{figure}
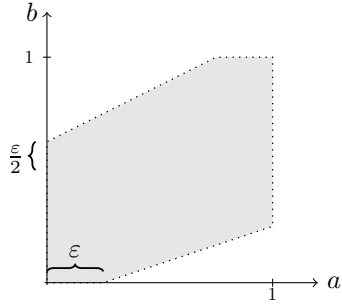

\medskip
In this note we aim to prove the following result.
\begin{theorem}\label{thm:theoremI}
Let $D$ be a bounded and Lipschitz domain contained in $\mathbb{R}^d$ and let $T>0$ be given. Let $\Omega=D\times (0,T)$ be the associated parabolic domain. Then there exists a positive number $\varepsilon\leq 1$ depending only on $D$ such that for each $p$ and each $s$ with $(s,\tfrac{1}{p})\in \mathcal{R}_{\varepsilon}$, a solution of the initial-boundary value problem
\begin{equation}
\left\{
  \begin{array}{ll}

\medskip
    \frac{\partial u}{\partial t}=\Delta u, & \hbox{in $\Omega$} \\
\medskip
    u(x,t)=f(x), & \hbox{for $(x,t)\in \partial D\times (0,T)$} \\\tag{$P$}
\medskip
    u(x,0)=g(x), & \hbox{for $x\in D$}
  \end{array}
\right.
\end{equation}
belongs to the parabolic Besov space $\mathbb{B}^{\alpha}_{\tau}(\Omega)$ with $0<\alpha<\min\bigl\{d\tfrac{p-1}{p},(s+\tfrac{1}{p})\tfrac{d}{d-1}\bigr\}$ and $\tfrac{1}{\tau}=\tfrac{\alpha}{d}+\tfrac{1}{p}$ provided that $f\in B^s_p(\partial D)$ and $g\in B^{s+\tfrac{1}{p}}_p(D)$.
\end{theorem}

\medskip
We shall precisely introduce in the next section all the Besov spaces involved in our main result, the elliptic Besov space $B^s_p(D)$ for a bounded Lipschitz domain in $\mathbb{R}^d$, the boundary elliptic Besov space $B^s_p(\partial D)$ and the parabolic Besov space $\mathbb{B}^s_p(\Omega)$ where $\Omega=D\times (0,T)$.

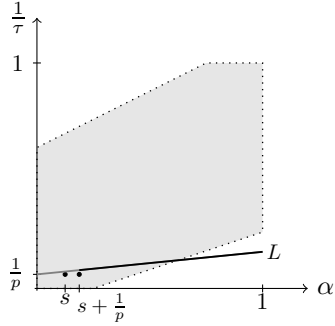
\begin{figure}[h!]
\begin{center}
\begin{tikzpicture}[scale=3]
\draw[line width=.5pt, dotted, fill=lightgray!40]
(0,0)--(1/4,0)--(1,1/4)--(1,1)--(3/4,1)--(0,5/8)--(0,0);
\draw[->] (-0.01,0)  -- (1.2,0) node[right] {\scalebox{1}{$\displaystyle \alpha$}};
\draw[->] (0,-0.01) -- (0,1.2) node[left] {\scalebox{1}{$\tfrac{1}{\tau}$}};
\draw[shift={(0,1)}] (.5pt,0pt) -- (-.5pt,0pt) node[left] {\scalebox{.9}{$1$}};
\draw[shift={(1,0)}] (0pt,-0.5pt)--(0pt,0.5pt) node[below] {\scalebox{.9}{$1$}};
\coordinate (s) at (1/8,1/16);
\filldraw [black]
     (s) circle (.25pt) node[above, black] {\scalebox{.8}{$$}};

\draw[shift={(0,1/16)}] (.5pt,0pt) -- (-.5pt,0pt) node[left] {\scalebox{.9}{$\tfrac{1}{p}$}};
\draw[shift={(1/8,0)}] (0pt,-0.5pt)--(0pt,0.5pt) node[below] {\scalebox{.9}{$s$}};

\coordinate (s+) at (1/8+1/16,1/16);
\filldraw [black]
     (s+) circle (.25pt) node[above right, black] {$$};
\coordinate (A) at (1/8+1/16+1/10,0.35);
\draw[shift={(1/8+1/16,0)}] (0pt,-0.5pt)--(0pt,0.5pt) node[below=of A] {\scalebox{.8}{$s+\tfrac{1}{p}$}};
\draw[line width=.7pt, gray] (0,1/16)--(1,1/16+1/10); 
\draw[line width=.7pt, black] (1/8+1/16,1/16+3/160)--(1,1/16+1/10);
\node at (1.05,1/16+1/10) {\scalebox{.8}{$L$}};
\end{tikzpicture}
\end{center}
\caption{The black segment in the line $L$ shows the improved regularity with data associated to $(s,\tfrac{1}{p})\in\mathcal{R}_{\varepsilon}$ when $d=10$.}\label{fig:improvedline}
\end{figure}

\medskip
The result becomes relevant when $\alpha$ can be taken larger than $s+\tfrac{1}{p}$, see Figure~\ref{fig:improvedline}. The arguments given in \cite{DaDeV97} regarding the improvement of the regularity parameter $\alpha$ as a tool to improve the rate of convergence of nonlinear approximation methods for elliptic problems, extend after Theorem~\ref{thm:theoremI}, to diffusion problems.

\section{Besov spaces}

This section is devoted to give a brief description of the Besov spaces involved in Theorem \ref{thm:theoremI}. Even when several approaches are possible, for the sake of simplicity we chose the interpolation one.

\medskip
The complete scale of Besov spaces $B^s_{p,q}(\mathbb{R}^d)$ indexed by the parameters $s$, $p$, $q$ in the space $\mathbb{R}^d$ is well known and several equivalent versions can be found in the classical literature such as Peetre's book \cite{Peetre76} and some modern approaches in the book by Y.~Meyer \cite{Meyer92} just to mention two standard references. We shall only deal with the case $p=q$ and we shall write $B^s_p$ instead of $B^s_{p,p}$. Given an open subset $D$ in $\mathbb{R}^d$ one can define $B^s_p(D)$ as the space of all the restrictions to $D$ of the functions in $B^s_p(\mathbb{R}^d)$. A second way to define $B^s_p(D)$ is provided by the real interpolation between Lebesgue and Sobolev spaces. Precisely, for $0<s<1$ and $1\leq p\leq \infty$, $B^s_p(D)=[L_p(D),W^1_p(D)]_{s,p}$ the $s$-interpolated between $L_p(D)$ and $W^1_p(D)$. When $D$ is a bounded Lipschitz domain in $\mathbb{R}^d$, both approaches coincide. The initial condition $g$ in our problem (P) belongs to $B^{s+\tfrac{1}{p}}_p(D)$ in the above described sense.

\medskip
The second Besov space is involved in the boundary condition $f$ and has to be described by using the local parametrization of the boundary $\partial D$ of $D$. After the standard localization arguments the problem of defining $B^s_p(\partial D)$, reduces to define the corresponding Besov class on the graph $G$ of a Lipschitz function $\phi$ with domain in $\mathbb{R}^{d-1}$. This is done by saying that $f\in B^s_p(G)$ when $f(x,\phi(x))\in B^s_p(\mathbb{R}^{d-1})$ for $0<s\leq 1$ and $p>0$.

\medskip
Let us now introduce through interpolation the parabolic Besov spaces involved in the statement of Theorem~\ref{thm:theoremI}. We would like to point out that Besov scales have been considered in very general settings such as spaces of homogeneous type, see \cite{HanSaw94} for example. The approach there is of Littlewood-Paley type. See also \cite{Schme77}, \cite{SchTri76} and \cite{BIN79} for more literature on the subject.
For $1\leq p\leq\infty$ the anisotropic Sobolev space $W^{2,1}_{p}(\Omega)$ is defined by the norm
\begin{equation*}
\norm{v}_{W^{2,1}_{p}(\Omega)}= \norm{v}_{L_p(\Omega)}
              + \sum_{i=1}^{d}\norm{\frac{\partial v}{\partial x_i}}_{L_p(\Omega)}
             + \sum_{i=1}^{d} \sum_{j=1}^{d}\norm{\frac{\partial^2 v}{\partial x_i\partial x_j}}_{L_p(\Omega)}
              + \norm{\frac{\partial v}{\partial t}}_{L_p(\Omega)}.
\end{equation*}
For $0<\alpha<2$ we define
\begin{equation*}
B^{\alpha,\tfrac{\alpha}{2}}_{p}(\Omega)=[L_p(\Omega),W^{2,1}_{p}(\Omega)]_{\tfrac{\alpha}{2},p},
\end{equation*}
the $\tfrac{\alpha}{2}$-real interpolated space between $W^{0,0}_{p}(\Omega)=L_p(\Omega)$ and $W^{2,1}_{p}(\Omega)$.
For simplicity we introduce the notation ${\mathbb{B}}^{\alpha}_{p}(\Omega)$ for the space $B^{\alpha,\tfrac{\alpha}{2}}_{p}(\Omega)$.

\section{Proof of Theorem~\ref{thm:theoremI}}

Let us start by stating two known regularity results for solutions of the heat equation. The first one, due to Jakab and Mitrea \cite{JaMit06}, gives the \textit{preservation} of Besov regularity from boundary and source data in the solution of the heat equation with homogeneous initial data. The second, contained in \cite{AGconstapp}, provides the \textit{improvement} of Besov regularity for the case of vanishing source.

\medskip
In the following statements, $D$ is a bounded Lipschitz domain in $\mathbb{R}^d$ ($d\geq 2$), $T>0$ and $\Omega=D\times (0,T)$.

\medskip
\textbf{Theorem JM} (Theorem (1.1) in \cite{JaMit06}).
There exists $\varepsilon>0$ depending only on the Lipschitz character of $D$ such that for every $(s,\tfrac{1}{p})\in\mathcal{R}_{\varepsilon}$ the solution of
\begin{equation}
\left\{
  \begin{array}{ll}
\medskip
    (\frac{\partial }{\partial t}-\Delta)v=\varphi, & \hbox{in $\Omega$} \\
\medskip
    v=\psi, & \hbox{in $\partial D\times (0,T)$} \\ \tag{$P_1$}
\medskip
    v(x,0)=0, & \hbox{for $x\in D$}
  \end{array}
\right.
\end{equation}
belongs to $L_p((0,T);B^{s+\tfrac{1}{p}}_{p}(D))$ provided that $\psi\in B^s_p(\partial D)$ and $\varphi\in B^{s+\tfrac{1}{p}-2}_{p}(D)$.

\medskip
\textbf{Theorem AG} (Theorem 1 and Theorem 2 in \cite{AGconstapp}).
For $1<p<\infty$, $\lambda>0$ and $u$ a solution of
\begin{equation*}
\frac{\partial u}{\partial t}=\Delta u \textrm{\quad in } \Omega
\end{equation*}
with $u\in L_p((0,T);B^{\lambda}_{p}(D))$, we have that $u\in\mathbb{B}^{\alpha}_{\tau}(\Omega)$ with $\tfrac{1}{\tau}=\tfrac{\alpha}{d}+\tfrac{1}{p}$ and
$0<\alpha<\min\{d\tfrac{p-1}{p},\lambda\tfrac{d}{d-1}\}$.

\medskip
Let $u$ be a solution of $(P)$ then $v=u-g$ is a solution of $(P_1)$ with $\varphi=-\Delta g$ and $\psi=f-g$. Now since $g\in B^{s+\tfrac{1}{p}}_p(D)$ we have that $\varphi=-\Delta g \in B^{s+\tfrac{1}{p}-2}_{p}(D)$. On the other hand from Theorem~3.1 in \cite{JeKe95} we know that the trace of $g$ on $\partial D$ belongs to $B^s_p(\partial D)$ since $g\in B^{s+\tfrac{1}{p}}_{p}(D)$. Hence $\psi\in B^s_p(\partial D)$. We are in position to apply Theorem~JM to obtain that $v\in L_p((0,T);B^{s+\tfrac{1}{p}}_{p}(D))$.

\medskip
Since $g$ belongs to $B^{s+\tfrac{1}{p}}_{p}(D)$ and $T$ is finite we have that $g\in L_p((0,T);B^{s+\tfrac{1}{p}}_{p}(D))$. Next we apply Theorem~AG with $\lambda=s+\tfrac{1}{p}$ and we get that $u\in \mathbb{B}^\alpha_\tau(\Omega)$ with $0<\alpha<\min\{d\tfrac{p-1}{p},\lambda\tfrac{d}{d-1}\}$ and $\tfrac{1}{\tau}=\tfrac{\alpha}{d}+\tfrac{1}{p}$ as desired.


\def\cprime{$'$} \def\cprime{$'$} \def\cprime{$'$}
\providecommand{\bysame}{\leavevmode\hbox to3em{\hrulefill}\thinspace}
\providecommand{\MR}{\relax\ifhmode\unskip\space\fi MR }
\providecommand{\MRhref}[2]{%
  \href{http://www.ams.org/mathscinet-getitem?mr=#1}{#2}
}
\providecommand{\href}[2]{#2}


\bigskip

{\footnotesize
\textsc{Instituto de Matem\'atica Aplicada del Litoral (IMAL), CONICET-UNL}


\medskip
\textmd{G\"{u}emes 3450}

\textmd{S3000GLN Santa Fe}

 \textmd{Argentina.}
\bigskip

\end{document}